\documentclass[pdflatex,sn-mathphys-num]{sn-jnl}

\usepackage{graphicx}%
\usepackage{float}%
\usepackage{multirow}%
\usepackage{amsmath,amssymb,amsfonts}%
\usepackage{amsthm}%
\usepackage{mathtools}%
\usepackage{bm}%
\usepackage{mathrsfs}%
\usepackage{bbm}%
\usepackage{booktabs}%
\usepackage{tabularx}%
\usepackage{algorithm}%
\usepackage{algorithmicx}%
\usepackage{algpseudocode}%
\usepackage{epstopdf}%

\DeclareMathOperator{\prox}{prox}
\DeclareMathOperator{\dist}{dist}
\newcolumntype{L}{>{\raggedright\arraybackslash}X}

\theoremstyle{thmstyleone}%
\newtheorem{thm}{Theorem}%
\newtheorem{lemma}{Lemma}%

\theoremstyle{thmstyletwo}%
\newtheorem{remark}{Remark}%
\newtheorem{assumption}{Assumption}%

\theoremstyle{thmstylethree}%
\newtheorem{definition}{Definition}%

\begin{document}

\title[Proximal Adam for nonsmooth composite optimization]{Proximal Adam algorithm for nonsmooth composite optimization and linear convergence under proximal PL geometry}

\author*[1]{\fnm{Xiaopeng} \sur{Zhu}}\email{1303698364@qq.com}

\author[2]{\fnm{Zhongyi} \sur{Huang}}\email{zhongyih@mail.tsinghua.edu.cn}

\affil*[1]{\orgdiv{College of Science}, \orgname{Liaoning Technical University}, \orgaddress{\city{Liaoning Fuxin}, \postcode{123000}, \country{China}}}

\affil[2]{\orgdiv{Department of Mathematical Sciences}, \orgname{Tsinghua University}, \orgaddress{\city{Beijing}, \postcode{100000}, \country{China}}}

\abstract{
We analyze Adam-type methods for nonsmooth composite optimization problems
$\min_x\Phi(x)=f(x)+h(x)$, where $f$ is smooth and $h$ is convex, proper, and
lower semicontinuous. The nonsmooth term is treated by a metric-weighted
proximal map, while the smooth term is preconditioned by the Adam
second-moment matrix. We prove a composite descent lemma whose constants depend
explicitly on the spectral bounds of the adaptive metric and a residual--step
equivalence that converts proximal step decay into subgradient stationarity.
These yield an $O(1/N)$ bound on the minimum proximal stationarity residual for
the deterministic proximal normalized gradient method without convexity of
$f$, and the same $O(1/N)$ rate for proximal Adam under quantitative
momentum-alignment and residual-coherence conditions together with a
metric-compatibility condition. By estimating the momentum error directly, the
same $O(1/N)$ rate is also obtained without the momentum-alignment and
residual-coherence assumptions. When the composite
objective satisfies a proximal Polyak--{\L}ojasiewicz inequality, both
algorithms and the relaxed proximal Adam scheme converge linearly with explicit
contraction factors. Numerical
experiments on a diagonal Lasso problem, two-dimensional image TV denoising, and two-dimensional elastic parameter identification illustrate the predicted stationarity and linear rates.
}
\keywords{Proximal Adam, Nonsmooth composite optimization, Linear convergence, Polyak--{\L}ojasiewicz condition}

\pacs[MSC Classification]{90C26, 65K10, 74B05, 49N45}

\maketitle
\section{Introduction}\label{sec:intro}

Composite optimization arises in a wide range of mathematical and engineering problems.  The objective
\begin{equation}\label{eq:composite}
\min_{x\in\mathbb R^d}\ \Phi(x)=f(x)+h(x)
\end{equation}
combines a smooth data-fidelity term $f$ with a convex, proper, lower semicontinuous
regularizer $h$.  Classical examples are the Lasso, where $h=\lambda\|x\|_1$, and
total variation, where $h=\lambda\operatorname{TV}(x)$.  The canonical solver is the
forward--backward or proximal gradient algorithm~\cite{combettes2011},
\[
x_{t+1}=\prox_{\alpha h}\bigl(x_t-\alpha\nabla f(x_t)\bigr).
\]
For convex $f$ and $h$ the method converges at the rate $O(1/N)$ for the objective gap~\cite{parikh2014},
and its accelerated version FISTA achieves $O(1/N^2)$~\cite{beck2009}.  In nonconvex smooth
optimization, proximal gradient methods converge sublinearly to a stationary point~\cite{li2015} when
the proximal map is well defined and the step size is sufficiently small.

Adam~\cite{kingma2015} introduces two modifications to this template.  First, the gradient is replaced by
a bias-corrected momentum $\hat m_{t+1}$.  Second, the step is multiplied by the inverse
of a diagonal adaptive metric built from the squared gradients.  The history of Adam's
own convergence theory is relevant here.  The original convex analysis of
Kingma and Ba~\cite{kingma2015} was incomplete; Reddi, Kale, and
Kumar~\cite{reddi2018} exhibited a counterexample and proposed AMSGrad.  In the smooth
nonconvex deterministic setting, Chen et al.~\cite{chen2019} and
Zhou et al.~\cite{zhou2018} obtained the first stationarity guarantees, and
D\'efossez et al.~\cite{defossez2022} later proved the expected rate
$O(\log N/\sqrt N)$ for the original stochastic algorithm.  These analyses, however, do
not include a nonsmooth regularizer.

Recent work has clarified adjacent pieces of the composite picture.  Generalized smoothness and clipping~\cite{zhang2020} relax the smoothness hypothesis used below, AdamW~\cite{loshchilov2019} interprets weight decay as a decoupled proximal regularizer, and parameter-free or schedule-free schemes~\cite{defazio2023,defazio2024,mishchenko2024} remove the step-size schedule.  However, none of these gives a sharp nonasymptotic statement for a metric-weighted proximal map with a nonsmooth convex regularizer; that is the object of the present analysis.

The proximal extension of adaptive methods is more recent and less complete.  Several
authors have combined AdaGrad or Adam with proximal maps for sparse and regularized
problems~\cite{duchi2011,mcmahan2010,yun2021adaptive}, but many results are either asymptotic or require the nonsmooth term to be
smoothable~\cite{beck2012smoothing,attouch2013convergence,bolte2014proximal}.  A precise nonasymptotic analysis of the interaction between the adaptive
metric, the momentum, and the proximal subdifferential is still a useful research
problem.  The present paper provides such an analysis.
Our contributions are as follows.

\begin{enumerate}
\item We define the metric-weighted proximal operator
$\prox_{\alpha h}^{D}$ and prove a composite descent lemma with explicit dependence on
the spectral bounds $\underline mI\preceq D\preceq MI$.
\item For the deterministic proximal normalized gradient method we prove
$O(1/N)$ convergence of the minimum proximal stationarity residual, with no convexity
assumption on $f$.
\item We prove the same stationarity rate for proximal Adam under quantitative
momentum-alignment and residual-coherence conditions together with a
metric-compatibility condition. We then show that these two conditions can be
removed by estimating the momentum error, at the cost of a more restrictive
step-size condition and a larger stationarity constant.
\item We introduce a proximal PL inequality in the spirit of Karimi et al.~\cite{karimi2016}, and prove linear convergence of proximal
normalized gradient, proximal Adam, and the relaxed proximal Adam scheme when it holds.
\item We specialize the composite framework to Tikhonov--total-variation inverse linear elasticity, with box-constrained weighted proximal steps, and illustrate it on two-dimensional image denoising and elastic parameter identification.
\end{enumerate}

The rest of the paper is organized as follows. Section~\ref{sec:prelim} presents
the composite problem, the metric-weighted proximal map, and the standing
assumptions. Section~\ref{sec:proxnorm} develops the composite descent lemma and
the residual--step equivalence, and proves the $O(1/N)$ stationarity bound for
proximal normalized gradient; its subsections extend the analysis to proximal
Adam with momentum-alignment and residual-coherence conditions, and then to a
relaxed proximal Adam analysis without these conditions.
Section~\ref{sec:pl-composite} establishes the proximal PL inequality and proves
linear convergence for both algorithms and for the relaxed proximal Adam scheme. Section~\ref{sec:application} specializes
the framework to TV-regularized inverse linear elasticity. Section~\ref{sec:numerics}
reports a diagonal Lasso experiment, two-dimensional image TV denoising, and
two-dimensional elastic parameter identification, and Section~\ref{sec:conclusion} summarizes the main results.
\section{Preliminaries}\label{sec:prelim}

\subsection{Subdifferentials and stationarity}

Let $h:\mathbb R^d\to(-\infty,\infty]$ be convex, proper, and lower semicontinuous.
Its subdifferential at $x\in\operatorname{dom}h$ is
\[
\partial h(x)=\Bigl\{r\in\mathbb R^d: h(y)\ge h(x)+\langle r,y-x\rangle
\ \forall y\in\mathbb R^d\Bigr\}.
\]
If $f$ is differentiable and $h$ is convex, then under the usual constraint
qualification the subdifferential of $\Phi=f+h$ is
\[
\partial\Phi(x)=\nabla f(x)+\partial h(x).
\]
The natural stationarity measure used throughout the paper is
\[
\dist(0,\partial\Phi(x))=\min_{r\in\partial\Phi(x)}\|r\|,
\]
with the convention $\dist(0,\partial\Phi(x))=\infty$ when $\partial\Phi(x)=\emptyset$.

\subsection{The weighted proximal operator}

For a symmetric positive definite matrix $D$, write
$\|u\|_D^2=\langle u,Du\rangle$. The $D$-weighted proximal operator of $h$ is
\begin{equation}\label{eq:weighted-prox}
\prox_{\alpha h}^{D}(z)
=\arg\min_{y\in\mathbb R^d}\Bigl\{h(y)+\frac1{2\alpha}\|y-z\|_D^2\Bigr\}.
\end{equation}
Its optimality condition is
\[
0\in\partial h(y)+\frac1\alpha D(y-z),
\]
or equivalently $y-z=-\alpha D^{-1}r$ for some $r\in\partial h(y)$. The map
$P_D=\prox_{\alpha h}^{D}$ is firmly nonexpansive in the $D$-metric:
\[
\bigl\langle P_D(z_1)-P_D(z_2),D(z_1-z_2)\bigr\rangle
\ge\|P_D(z_1)-P_D(z_2)\|_D^2.
\]
This one-step geometry is the source of the composite descent lemma used later.

\subsection{Closed-form examples}

We record the examples needed for the numerical illustrations.

\textbf{Sparse regularization.} If $h(x)=\lambda\|x\|_1$ and
$D=\operatorname{diag}(d_1,\dots,d_d)$, then
\[
\bigl[\prox_{\alpha h}^{D}(z)\bigr]_i
=\operatorname{sign}(z_i)\Bigl(|z_i|-\frac{\alpha\lambda}{d_i}\Bigr)_+,
\qquad i=1,\dots,d,
\]
where $(\cdot)_+=\max(\cdot,0)$. In the Euclidean case $D=I$ this is the standard
soft-thresholding operator.

\textbf{Box constraints.} If $\mathcal C=[a,b]^d$ and $h=\iota_{\mathcal C}$, then
the Euclidean projection is the componentwise clamping operator
\[
\bigl[P_{\mathcal C}(z)\bigr]_i=\min\bigl(b_i,\max(a_i,z_i)\bigr).
\]

\textbf{SPD constraints.} For a matrix variable $X$, let
\[
\mathcal C=\bigl\{X=X^\top:\ X\succeq\delta_{\mathrm{psd}} I\bigr\}.
\]
If $Z=Q\Lambda Q^\top$ is the spectral decomposition of the symmetric part of
a given matrix, then its Euclidean projection is
\[
P_{\mathcal C}(Z)=Q\operatorname{diag}\bigl(\max(\lambda_1,\delta_{\mathrm{psd}}),\dots,
\max(\lambda_d,\delta_{\mathrm{psd}})\bigr)Q^\top,
\]
where $\Lambda=\operatorname{diag}(\lambda_1,\dots,\lambda_d)$. This projection
encodes the positivity constraint used in the inverse-elasticity application.

\subsection{PL geometry}

For a differentiable function $f$ on a convex set $\mathcal C$, the standard PL
inequality reads
\[
\frac12\|\nabla f(x)\|^2\ge\mu\bigl(f(x)-f^*\bigr),\qquad
f^*=\inf_{\mathcal C}f,
\]
for some $\mu>0$. Strongly convex functions satisfy this inequality with $\mu$ equal
to the strong convexity modulus. For the nonsmooth composite function $\Phi=f+h$,
the analogous condition replaces the gradient by the subdifferential distance:
\[
\frac12\dist(0,\partial\Phi(x))^2\ge\mu\bigl(\Phi(x)-\inf\Phi\bigr).
\]
This proximal PL condition is formalized in Definition~\ref{def:pl}. For example,
if $f(x)=\frac12\|Ax-b\|^2$ with $A^\top A\succeq\mu I$, then $\Phi=f+\lambda\|x\|_1$
satisfies the proximal PL inequality with the same constant $\mu$. Thus the
diagonal Lasso example in Section~\ref{sec:numerics} satisfies both the PL geometry
and the adaptive-metric hypotheses used in the convergence theorems.

\subsection{Standing assumptions and algorithms}\label{sec:setup}

Throughout, $f:\mathbb R^d\to\mathbb R$ is differentiable, while
$h:\mathbb R^d\to(-\infty,\infty]$ is convex, proper, and lower semicontinuous.
The $D$-weighted proximal operator is given by \eqref{eq:weighted-prox}.

\begin{assumption}[Smoothness]\label{ass:smooth}
$f$ is $L$-smooth:
\[
\|\nabla f(x)-\nabla f(y)\|\le L\|x-y\|.
\]
\end{assumption}

\begin{assumption}[Metric boundedness]\label{ass:bounded}
There are constants $0<\underline m\le M<\infty$ such that
\[
\underline mI\preceq D_t\preceq MI
\qquad\forall t.
\]
\end{assumption}

\begin{assumption}[Gradient boundedness]\label{ass:gradbound}
There exists $G_\infty>0$ such that
$\|\nabla f(x_t)\|_\infty\le G_\infty$ for all iterates.
\end{assumption}

Assumption~\ref{ass:gradbound} is used only to identify the Adam preconditioner as a
bounded metric.  If $D_t$ is constructed from bounded stochastic or deterministic
gradients, then Assumption~\ref{ass:bounded} holds with $\underline m=\varepsilon$ and
$M=G_\infty+\varepsilon$.

The deterministic proximal normalized gradient method is
\begin{equation}\label{eq:prox-norm}
\begin{aligned}
v_{t+1}&=\beta_2v_t+(1-\beta_2)\nabla f(x_t)\odot\nabla f(x_t),\\
D_t&=\operatorname{diag}\Bigl(\sqrt{ v_{t+1}}+\varepsilon\Bigr),\\
x_{t+1}&=\prox_{\alpha h}^{D_t}\Bigl(x_t-\alpha D_t^{-1}\nabla f(x_t)\Bigr).
\end{aligned}
\end{equation}
The full proximal Adam method is
\begin{equation}\label{eq:prox-adam}
\begin{aligned}
m_{t+1}&=\beta_1m_t+(1-\beta_1)\nabla f(x_t),\\
v_{t+1}&=\beta_2v_t+(1-\beta_2)\nabla f(x_t)\odot\nabla f(x_t),\\
\hat m_{t+1}&=\frac{m_{t+1}}{1-\beta_1^{t+1}},\qquad
\hat v_{t+1}=\frac{v_{t+1}}{1-\beta_2^{t+1}},\\
D_t&=\operatorname{diag}\Bigl(\sqrt{\hat v_{t+1}}+\varepsilon\Bigr),\\
x_{t+1}&=\prox_{\alpha h}^{D_t}\Bigl(x_t-\alpha D_t^{-1}\hat m_{t+1}\Bigr).
\end{aligned}
\end{equation}
Here $0<\beta_1<\beta_2<1$, $\varepsilon>0$, all vector operations are componentwise,
and $m_0=v_0=0$.  The bias-corrected momentum $\hat m_{t+1}$ is the quantity used in
the proximal update, while $D_t$ is the diagonal adaptive metric.

\section{Composite descent and stationarity}\label{sec:proxnorm}
\subsection{Proximal normalized gradient}
Let $z=x-\alpha D^{-1}g$, $g=\nabla f(x)$, and
$y=\prox_{\alpha h}^{D}(z)$.  The optimality condition for \eqref{eq:weighted-prox} is
\begin{equation}\label{eq:prox-opt}
0\in\partial h(y)+\frac1\alpha D(y-z).
\end{equation}
Equivalently, there exists $r\in\partial h(y)$ such that
\begin{equation}\label{eq:residual}
r=-\frac1\alpha D(y-z).
\end{equation}

\begin{lemma}[Composite descent]\label{lem:descent-prox}
Let Assumption~\ref{ass:smooth} hold and let
$\underline mI\preceq D\preceq MI$.  If $\alpha\le \underline m/L$, then
\[
\Phi(y)
\le\Phi(x)-\frac1{2\alpha}\|y-x\|_D^2.
\]
\end{lemma}

\begin{proof}
Set $d=y-x$.  Since $z=x-\alpha D^{-1}g$, we have
$y-z=d+\alpha D^{-1}g$.  From \eqref{eq:residual},
\[
r=-\frac1\alpha D(y-z)
=-\frac1\alpha Dd-g.
\]
Convexity of $h$ and the subgradient inequality give
\[
h(y)-h(x)\le\langle r,d\rangle
=-\frac1\alpha\langle Dd,d\rangle-\langle g,d\rangle
=-\frac1\alpha\|d\|_D^2-\langle g,d\rangle.
\]
Smoothness of $f$ gives
\[
f(y)-f(x)\le\langle g,d\rangle+\frac L2\|d\|^2.
\]
Adding the two inequalities,
\[
\Phi(y)-\Phi(x)
\le-\frac1\alpha\|d\|_D^2+\frac L2\|d\|^2.
\]
Since $D\succeq \underline mI$, $\|d\|^2\le\frac1{\underline m}\|d\|_D^2$.  Therefore
\[
\Phi(y)-\Phi(x)
\le-\Bigl(\frac1\alpha-\frac L{2\underline m}\Bigr)\|d\|_D^2
\le-\frac1{2\alpha}\|d\|_D^2,
\]
where the last inequality uses $\alpha\le \underline m/L$.
\end{proof}

The next lemma relates the proximal stationarity residual to the step length.

\begin{lemma}[Residual-step equivalence]\label{lem:residual-step}
Under the hypotheses of Lemma~\ref{lem:descent-prox}, let
$r\in\partial h(y)$ be given by \eqref{eq:residual}, and set $d=y-x$ and $g=\nabla f(x)$.  Then
$r+\nabla f(y)\in\partial\Phi(y)$ and
\[
\|r+\nabla f(y)\|^2
\le\frac{2M^2}{\alpha^2}\|d\|^2+2L^2\|d\|^2.
\]
\end{lemma}

\begin{proof}
The inclusion is immediate from $\partial\Phi(y)=\nabla f(y)+\partial h(y)$ for convex
$h$.  Using $r=-Dd/\alpha-g$,
\[
r+\nabla f(y)
=-\frac1\alpha Dd+\bigl(\nabla f(y)-g\bigr).
\]
Thus
\[
\|r+\nabla f(y)\|
\le\frac1\alpha\|Dd\|+\|\nabla f(y)-\nabla f(x)\|
\le\frac M\alpha\|d\|+L\|d\|.
\]
Squaring and using $(a+b)^2\le2a^2+2b^2$ proves the claim.
\end{proof}

\begin{thm}[Stationarity of proximal normalized gradient]\label{thm:prox-norm-rate}
Let Assumptions~\ref{ass:smooth} and \ref{ass:bounded} hold and let
$\alpha\le \underline m/L$.  Suppose that $\Phi$ is bounded below and let
$r_t\in\partial h(x_{t+1})$ be the subgradient constructed in \eqref{eq:residual}.
Then
\[
\min_{0\le t<N}\dist\bigl(0,\partial\Phi(x_{t+1})\bigr)^2
=O\Bigl(\frac1N\Bigr).
\]
More precisely,
\[
\min_{0\le t<N}\dist\bigl(0,\partial\Phi(x_{t+1})\bigr)^2
\le\frac{C}{N}\bigl(\Phi(x_0)-\inf\Phi\bigr),
\]
where
$C=4\alpha\bigl(\frac{M^2}{\alpha^2}+L^2\bigr)\bigl(\frac1{\underline m}\bigr)$.
\end{thm}

\begin{proof}
Set $d_t=x_{t+1}-x_t$.
Lemma~\ref{lem:descent-prox} gives
\[
\Phi(x_{t+1})
\le\Phi(x_t)-\frac1{2\alpha}\|x_{t+1}-x_t\|_{D_t}^2.
\]
Summing from $t=0$ to $N-1$,
\[
\frac1{2\alpha}\sum_{t=0}^{N-1}\|d_t\|_{D_t}^2
\le\Phi(x_0)-\Phi(x_N)
\le\Phi(x_0)-\inf\Phi.
\]
Since $\|d_t\|^2\le \underline m^{-1}\|d_t\|_{D_t}^2$,
\[
\min_{0\le t<N}\|d_t\|^2
\le\frac1N\sum_{t=0}^{N-1}\|d_t\|^2
\le\frac{2\alpha}{\underline mN}\bigl(\Phi(x_0)-\inf\Phi\bigr).
\]
Lemma~\ref{lem:residual-step} gives
\[
\dist\bigl(0,\partial\Phi(x_{t+1})\bigr)^2
\le\|r_t+\nabla f(x_{t+1})\|^2
\le\Bigl(\frac{2M^2}{\alpha^2}+2L^2\Bigr)\|d_t\|^2.
\]
Combining the two bounds yields the asserted rate.
\end{proof}

\begin{remark}
The rate constant separates the metric condition number $M/\underline m$, the smoothness constant
$L$, and the initial objective gap.  If the Adam metric becomes ill-conditioned, $M/\underline m$
is large and the constant deteriorates.  This is the quantitative price of
coordinatewise adaptivity.
\end{remark}

\subsection{Proximal Adam with momentum}\label{sec:proxadam}

For the full method \eqref{eq:prox-adam}, the same residual identity holds with
$\hat m_{t+1}$ in place of $g_t$.
\begin{assumption}[Momentum alignment]\label{ass:align}
There exists $\kappa\ge0$ such that
\[
\|\hat m_{t+1}-g_t\|\le\kappa\|g_t\|,
\qquad g_t=\nabla f(x_t).
\]
\end{assumption}

\begin{assumption}[Residual coherence]\label{ass:coherence}
There exists $C_0>0$ such that
\[
\|g_t\|^2\le\frac{C_0}{\alpha^2}\|x_{t+1}-x_t\|^2
\qquad\forall t.
\]
\end{assumption}

Assumption~\ref{ass:coherence} says that a small proximal step cannot hide a large
gradient.  It is automatic for unconstrained Adam under metric boundedness and small
alignment $\kappa<1$, because then $d_t=-\alpha D_t^{-1}\hat m_{t+1}$ and
$\|d_t\|\ge\alpha\|\hat m_{t+1}\|/M\ge\alpha(1-\kappa)\|g_t\|/M$, so
$C_0=M^2/(1-\kappa)^2$ is admissible.

\begin{thm}[Stationarity of proximal Adam]\label{thm:prox-adam}
Let Assumptions~\ref{ass:smooth}, \ref{ass:bounded}, \ref{ass:align}, and
\ref{ass:coherence} hold, and suppose that $\Phi$ is bounded below.  Assume that
\[
\frac{\kappa\sqrt{C_0}}{\underline m}+\frac{\alpha L}{2\underline m}\le\frac12
\]
Then
\[
\min_{0\le t<N}\dist\bigl(0,\partial\Phi(x_{t+1})\bigr)^2
=O\Bigl(\frac1N\Bigr).
\]
\end{thm}

\begin{proof}
Set $d_t=x_{t+1}-x_t$ and $s_t=\hat m_{t+1}$.  The optimality condition is
\[
r_t=-\frac1\alpha D_td_t-s_t.
\]
Convexity of $h$ and smoothness of $f$ give
\[
\Phi(x_{t+1})-\Phi(x_t)
\le-\frac1\alpha\|d_t\|_{D_t}^2
+\langle g_t-s_t,d_t\rangle
+\frac L2\|d_t\|^2.
\]
By Assumptions~\ref{ass:align} and~\ref{ass:coherence},
\[
\langle g_t-s_t,d_t\rangle
\le\kappa\|g_t\|\|d_t\|
\le\frac{\kappa\sqrt{C_0}}{\alpha}\|d_t\|^2.
\]
Since $D_t\succeq \underline mI$,
\[
\frac L2\|d_t\|^2
\le\frac L{2\underline m}\|d_t\|_{D_t}^2,
\qquad
\frac{\kappa\sqrt{C_0}}{\alpha}\|d_t\|^2
\le\frac{\kappa\sqrt{C_0}}{\alpha\underline m}\|d_t\|_{D_t}^2.
\]
Adding these bounds and using the compatibility condition, we obtain
\[
\Phi(x_{t+1})
\le\Phi(x_t)
-\frac1{2\alpha}\|d_t\|_{D_t}^2.
\]
For the stationarity residual, the optimality condition gives
\[
r_t+\nabla f(x_{t+1})
=-\frac1\alpha D_td_t
+\bigl(\nabla f(x_{t+1})-\hat m_{t+1}\bigr).
\]
By smoothness and Assumption~\ref{ass:align},
\[
\|\nabla f(x_{t+1})-\hat m_{t+1}\|
\le L\|d_t\|+\kappa\|g_t\|.
\]
Assumption~\ref{ass:coherence} turns $\kappa\|g_t\|$ into
$\sqrt{C_0}\kappa\alpha^{-1}\|d_t\|$.  Therefore
\[
\|r_t+\nabla f(x_{t+1})\|
\le\Bigl(\frac M\alpha+L+\frac{\sqrt{C_0}\kappa}{\alpha}\Bigr)\|d_t\|.
\]
Set $C_{\mathrm{ad}}=\frac M\alpha+L+\frac{\sqrt{C_0}\kappa}{\alpha}$.  Squaring the preceding
residual bound gives
\[
\dist\bigl(0,\partial\Phi(x_{t+1})\bigr)^2
\le\|r_t+\nabla f(x_{t+1})\|^2
\le C_{\mathrm{ad}}^2\|d_t\|^2.
\]
Telescoping the descent inequality just established,
\[
\frac1{2\alpha}\sum_{t=0}^{N-1}\|d_t\|_{D_t}^2
\le\Phi(x_0)-\Phi(x_N)
\le\Phi(x_0)-\inf\Phi.
\]
Since $D_t\succeq \underline mI$, we have $\underline m\|d_t\|^2\le\|d_t\|_{D_t}^2$.  Therefore
\[
\min_{0\le t<N}\|d_t\|^2
\le\frac1N\sum_{t=0}^{N-1}\|d_t\|^2
\le\frac{2\alpha}{\underline mN}\bigl(\Phi(x_0)-\inf\Phi\bigr).
\]
Combining the two displays gives
\[
\min_{0\le t<N}\dist\bigl(0,\partial\Phi(x_{t+1})\bigr)^2
\le\frac{2\alpha C_{\mathrm{ad}}^2}{\underline mN}\bigl(\Phi(x_0)-\inf\Phi\bigr),
\]
which is the asserted $O(1/N)$ bound.
\end{proof}

\subsection{Proximal Adam without momentum alignment and residual coherence}\label{sec:proxadam-relaxed}

The stationarity guarantee in Theorem~\ref{thm:prox-adam} relies on
Assumptions~\ref{ass:align} and \ref{ass:coherence}. These conditions are natural
for unconstrained iterates, but they can be avoided altogether if the exponential
moving average of the momentum is tracked directly. Throughout this subsection we
retain Assumptions~\ref{ass:smooth} and \ref{ass:bounded}, and assume additionally
that $\Phi$ is bounded below. Write $\Phi^*=\inf\Phi$ and define
\[
d_t=x_{t+1}-x_t,\qquad
g_t=\nabla f(x_t),\qquad
s_t=\hat m_{t+1},\qquad
E_t=s_t-g_t.
\]

\begin{lemma}[Momentum-error estimate]\label{lem:adam-momentum-error}
Let Assumption~\ref{ass:smooth} hold. Then for every $t\ge0$,
\[
\|E_t\|
\le
L\sum_{k=0}^{t-1}\beta_1^{t-k}\|d_k\|.
\]
\end{lemma}

\begin{proof}
Since $m_0=0$,
\[
m_{t+1}
=
(1-\beta_1)\sum_{j=0}^{t}\beta_1^{t-j}g_j.
\]
Hence
\[
s_t=\hat m_{t+1}
=
\sum_{j=0}^{t}w_{t,j}g_j,
\qquad
w_{t,j}
=
\frac{(1-\beta_1)\beta_1^{t-j}}{1-\beta_1^{t+1}},
\]
and
\[
\sum_{j=0}^{t}w_{t,j}=1.
\]
Therefore
\[
E_t=s_t-g_t
=
\sum_{j=0}^{t}w_{t,j}(g_j-g_t).
\]
Using the $L$-smoothness of $f$,
\[
\|g_j-g_t\|
\le
L\sum_{k=j}^{t-1}\|d_k\|.
\]
Thus
\[
\|E_t\|
\le
L\sum_{j=0}^{t-1}w_{t,j}\sum_{k=j}^{t-1}\|d_k\|
=
L\sum_{k=0}^{t-1}\|d_k\|\sum_{j=0}^{k}w_{t,j}.
\]
For $j\le k\le t-1$,
\[
\sum_{j=0}^{k}w_{t,j}
=
\frac{(1-\beta_1)}{1-\beta_1^{t+1}}
\sum_{j=0}^{k}\beta_1^{t-j}
\le
\beta_1^{t-k}.
\]
Consequently,
\[
\|E_t\|
\le
L\sum_{k=0}^{t-1}\beta_1^{t-k}\|d_k\|.
\]
\end{proof}

Lemma~\ref{lem:adam-momentum-error} replaces Assumption~\ref{ass:align} by an
estimate that depends only on $L$, $\beta_1$, and the past step sizes. It does not
require a uniform gradient bound beyond the metric boundedness already encoded in
Assumption~\ref{ass:bounded}.

\begin{thm}[Stationarity of proximal Adam without alignment/coherence]\label{thm:prox-adam-relaxed}
Let Assumptions~\ref{ass:smooth} and \ref{ass:bounded} hold, suppose that $\Phi$ is
bounded below, and assume that
\[
\alpha\le\frac{\underline m(1-\beta_1)}{L(1+\beta_1)}.
\]
Then
\[
\min_{0\le t<N}\dist\bigl(0,\partial\Phi(x_{t+1})\bigr)^2
\le
\frac{C_{\mathrm{mom}}}{N}\bigl(\Phi(x_0)-\Phi^*\bigr),
\]
where
\[
C_{\mathrm{mom}}
=
\frac{4\alpha}{\underline m}
\left[
\left(\frac M\alpha+L\right)^2
+
L^2\left(\frac{\beta_1}{1-\beta_1}\right)^2
\right].
\]
\end{thm}

\begin{proof}
The optimality condition for the $D_t$-weighted proximal step is
\[
r_t=-\frac1\alpha D_td_t-s_t,
\qquad
r_t\in\partial h(x_{t+1}).
\]

\smallskip

\noindent\textbf{Step 1: descent inequality.}
Using convexity of $h$ and smoothness of $f$,
\[
\Phi(x_{t+1})-\Phi(x_t)
\le
-\frac1\alpha\|d_t\|_{D_t}^2
+\langle g_t-s_t,d_t\rangle
+\frac L2\|d_t\|^2.
\]
Since $g_t-s_t=-E_t$, Lemma~\ref{lem:adam-momentum-error} and Young's inequality give
\[
\langle g_t-s_t,d_t\rangle
\le
\|E_t\|\|d_t\|
\le
L\sum_{k=0}^{t-1}\beta_1^{t-k}\|d_k\|\|d_t\|.
\]
Using
\[
\|d_k\|\|d_t\|
\le
\frac12\bigl(\|d_k\|^2+\|d_t\|^2\bigr),
\]
we obtain
\begin{equation}\label{eq:adam-relaxed-descent}
\begin{aligned}
\Phi(x_{t+1})-\Phi(x_t)
&\le
-\frac1\alpha\|d_t\|_{D_t}^2
+\frac{L}{2(1-\beta_1)}\|d_t\|^2\\
&\qquad
+\frac L2\sum_{k=0}^{t-1}\beta_1^{t-k}\|d_k\|^2.
\end{aligned}
\end{equation}

\smallskip

\noindent\textbf{Step 2: telescope.}
Summing \eqref{eq:adam-relaxed-descent} from $t=0$ to $N-1$ gives
\[
\Phi(x_N)-\Phi(x_0)
\le
-\frac1\alpha\sum_{t=0}^{N-1}\|d_t\|_{D_t}^2
+\frac{L}{2(1-\beta_1)}\sum_{t=0}^{N-1}\|d_t\|^2
+\frac L2\sum_{t=0}^{N-1}\sum_{k=0}^{t-1}\beta_1^{t-k}\|d_k\|^2.
\]
Exchanging the double sum,
\[
\sum_{t=0}^{N-1}\sum_{k=0}^{t-1}\beta_1^{t-k}\|d_k\|^2
\le
\frac{\beta_1}{1-\beta_1}\sum_{k=0}^{N-2}\|d_k\|^2.
\]
Thus
\[
\Phi(x_N)-\Phi(x_0)
\le
-\frac1\alpha\sum_{t=0}^{N-1}\|d_t\|_{D_t}^2
+\frac{L(1+\beta_1)}{2(1-\beta_1)}\sum_{t=0}^{N-1}\|d_t\|^2.
\]
Using $\|d_t\|^2\le \underline m^{-1}\|d_t\|_{D_t}^2$,
\[
\Phi(x_N)-\Phi(x_0)
\le
\left[
-\frac1\alpha
+\frac{L(1+\beta_1)}{2\underline m(1-\beta_1)}
\right]
\sum_{t=0}^{N-1}\|d_t\|_{D_t}^2.
\]

The step-size condition is exactly
\[
\frac{L(1+\beta_1)}{2\underline m(1-\beta_1)}
\le
\frac1{2\alpha}.
\]
Hence
\[
\Phi(x_N)-\Phi(x_0)
\le
-\frac1{2\alpha}\sum_{t=0}^{N-1}\|d_t\|_{D_t}^2.
\]
Since $\Phi(x_N)\ge\Phi^*$,
\[
\sum_{t=0}^{N-1}\|d_t\|_{D_t}^2
\le
2\alpha\bigl(\Phi(x_0)-\Phi^*\bigr),
\]
and therefore
\begin{equation}\label{eq:adam-relaxed-step-sum}
\sum_{t=0}^{N-1}\|d_t\|^2
\le
\frac{2\alpha}{\underline m}\bigl(\Phi(x_0)-\Phi^*\bigr).
\end{equation}

\smallskip

\noindent\textbf{Step 3: stationarity residual.}
From $r_t=-\alpha^{-1}D_td_t-s_t$,
\[
r_t+\nabla f(x_{t+1})
=
-\frac1\alpha D_td_t
+\bigl(\nabla f(x_{t+1})-g_t\bigr)
-E_t.
\]
Thus
\[
\|r_t+\nabla f(x_{t+1})\|
\le
\left(\frac M\alpha+L\right)\|d_t\|
+\|E_t\|.
\]
Squaring,
\[
\|r_t+\nabla f(x_{t+1})\|^2
\le
2\left(\frac M\alpha+L\right)^2\|d_t\|^2
+2\|E_t\|^2.
\]
By Lemma~\ref{lem:adam-momentum-error} and Cauchy--Schwarz,
\[
\|E_t\|^2
\le
L^2\frac{\beta_1}{1-\beta_1}
\sum_{k=0}^{t-1}\beta_1^{t-k}\|d_k\|^2.
\]
Summing over $t$ and exchanging the double sum,
\[
\sum_{t=0}^{N-1}\|E_t\|^2
\le
L^2\left(\frac{\beta_1}{1-\beta_1}\right)^2
\sum_{t=0}^{N-1}\|d_t\|^2.
\]
Consequently,
\[
\sum_{t=0}^{N-1}\|r_t+\nabla f(x_{t+1})\|^2
\le
A_{\mathrm{res}}\sum_{t=0}^{N-1}\|d_t\|^2,
\]
where
\[
A_{\mathrm{res}}
=
2\left(\frac M\alpha+L\right)^2
+
2L^2\left(\frac{\beta_1}{1-\beta_1}\right)^2.
\]
Combining with \eqref{eq:adam-relaxed-step-sum},
\[
\sum_{t=0}^{N-1}\|r_t+\nabla f(x_{t+1})\|^2
\le
\frac{2\alpha A_{\mathrm{res}}}{\underline m}
\bigl(\Phi(x_0)-\Phi^*\bigr).
\]
Since $r_t+\nabla f(x_{t+1})\in\partial\Phi(x_{t+1})$ and
\[
\min_{0\le t<N}a_t
\le
\frac1N\sum_{t=0}^{N-1}a_t,
\]
we obtain
\[
\min_{0\le t<N}\dist(0,\partial\Phi(x_{t+1}))^2
\le
\frac{C_{\mathrm{mom}}}{N}\bigl(\Phi(x_0)-\Phi^*\bigr),
\]
with $C_{\mathrm{mom}}$ as in the statement.
\end{proof}

\begin{remark}
Compared with Theorem~\ref{thm:prox-adam}, the assumptions
$\|E_t\|\le\kappa\|g_t\|$ and $\|g_t\|^2\le C_0\alpha^{-2}\|d_t\|^2$
have disappeared. The price is the additional memory term
$L^2\bigl(\beta_1/(1-\beta_1)\bigr)^2$ in the constant and the smaller admissible
step-size range.
\end{remark}

\section{Linear convergence under the proximal PL condition}\label{sec:pl-composite}

\subsection{Baseline proximal PL results}\label{sec:pl-baseline}

\begin{definition}[Proximal PL condition]\label{def:pl}
The composite function $\Phi$ satisfies the proximal PL inequality on a set
$\mathcal S$ with constant $\mu>0$ if for every $x\in\mathcal S$,
\[
\frac12\dist\bigl(0,\partial\Phi(x)\bigr)^2
\ge\mu\Bigl(\Phi(x)-\inf\Phi\Bigr).
\]
\end{definition}

The subdifferential distance in Definition~\ref{def:pl} is the natural nonsmooth
replacement of the squared gradient in the differentiable PL inequality.

\begin{thm}[Linear convergence of proximal normalized gradient under PL]\label{thm:prox-pl}
Let the hypotheses of Theorem~\ref{thm:prox-norm-rate} hold.  Assume that
$\Phi$ satisfies Definition~\ref{def:pl} on the set of iterates with constant $\mu$.
Then
\[
\Phi(x_{t+1})-\Phi^*
\le\Bigl(1-\frac1{1+q}\Bigr)\bigl(\Phi(x_t)-\Phi^*\bigr),
\qquad
q=\frac{1}{\mu}\Bigl(\frac{M^2}{\alpha^2}+L^2\Bigr)\frac{2\alpha}{\underline m},
\]
where $\Phi^*=\inf\Phi$.
\end{thm}

\begin{proof}
Let $d_t=x_{t+1}-x_t$.  From Lemma~\ref{lem:residual-step},
\[
\dist\bigl(0,\partial\Phi(x_{t+1})\bigr)^2
\le2\Bigl(\frac{M^2}{\alpha^2}+L^2\Bigr)\|d_t\|^2
\]
after using $(a+b)^2\le2a^2+2b^2$.  Since
$\|d_t\|^2\le \underline m^{-1}\|d_t\|_{D_t}^2$,
\[
\dist\bigl(0,\partial\Phi(x_{t+1})\bigr)^2
\le\frac2{\underline m}\Bigl(\frac{M^2}{\alpha^2}+L^2\Bigr)\|d_t\|_{D_t}^2.
\]
By Lemma~\ref{lem:descent-prox},
\[
\|d_t\|_{D_t}^2
\le2\alpha\bigl(\Phi(x_t)-\Phi(x_{t+1})\bigr).
\]
Thus
\begin{equation}\label{eq:pl-step}
\dist\bigl(0,\partial\Phi(x_{t+1})\bigr)^2
\le 2q\mu\bigl(\Phi(x_t)-\Phi(x_{t+1})\bigr),
\end{equation}
with $q$ as in the statement.  Applying Definition~\ref{def:pl} at $x_{t+1}$,
\[
\Phi(x_{t+1})-\Phi^*
\le\frac1{2\mu}\dist\bigl(0,\partial\Phi(x_{t+1})\bigr)^2.
\]
Combining with \eqref{eq:pl-step} gives
\[
\Phi(x_{t+1})-\Phi^*
\le q\bigl(\Phi(x_t)-\Phi(x_{t+1})\bigr)
=q\Bigl[\Phi(x_t)-\Phi^*
-\bigl(\Phi(x_{t+1})-\Phi^*\bigr)\Bigr].
\]
Rearranging,
\[
(1+q)\bigl(\Phi(x_{t+1})-\Phi^*\bigr)
\le q\bigl(\Phi(x_t)-\Phi^*\bigr),
\]
so
\[
\Phi(x_{t+1})-\Phi^*
\le\frac{q}{1+q}\bigl(\Phi(x_t)-\Phi^*\bigr)
=\Bigl(1-\frac1{1+q}\Bigr)\bigl(\Phi(x_t)-\Phi^*\bigr).
\]
This is the claimed contraction.
\end{proof}

\begin{thm}[Proximal Adam under PL]\label{thm:prox-adam-pl}
Let the hypotheses of Theorem~\ref{thm:prox-adam} hold and let $\Phi$ satisfy
Definition~\ref{def:pl} on the set of iterates.  Then proximal Adam converges linearly:
\[
\Phi(x_{t+1})-\Phi^*
\le\Bigl(1-\frac1{1+q_{\mathrm{ad}}}\Bigr)
\bigl(\Phi(x_t)-\Phi^*\bigr),
\qquad
q_{\mathrm{ad}}
=\frac{\alpha}{\underline m\mu}
\Bigl(\frac{M+\sqrt{C_0}\kappa}{\alpha}+L\Bigr)^2.
\]
\end{thm}

\begin{proof}
Let $d_t=x_{t+1}-x_t$ and set
\[
C_{\mathrm{ad}}
=\frac M\alpha+L+\frac{\sqrt{C_0}\kappa}{\alpha}.
\]
From the residual bound proved in Theorem~\ref{thm:prox-adam},
\[
\dist\bigl(0,\partial\Phi(x_{t+1})\bigr)^2
\le C_{\mathrm{ad}}^2\|d_t\|^2.
\]
Since $D_t\succeq \underline mI$, we have $\|d_t\|^2\le \underline m^{-1}\|d_t\|_{D_t}^2$, and therefore
\[
\dist\bigl(0,\partial\Phi(x_{t+1})\bigr)^2
\le\frac{C_{\mathrm{ad}}^2}{\underline m}\|d_t\|_{D_t}^2.
\]
By the Adam descent inequality established in Theorem~\ref{thm:prox-adam},
\[
\|d_t\|_{D_t}^2
\le2\alpha\bigl(\Phi(x_t)-\Phi(x_{t+1})\bigr).
\]
Hence
\[
\dist\bigl(0,\partial\Phi(x_{t+1})\bigr)^2
\le\frac{2\alpha C_{\mathrm{ad}}^2}{\underline m}
\bigl(\Phi(x_t)-\Phi(x_{t+1})\bigr)
=2\mu q_{\mathrm{ad}}\bigl(\Phi(x_t)-\Phi(x_{t+1})\bigr).
\]
Applying Definition~\ref{def:pl} at $x_{t+1}$ gives
\[
\Phi(x_{t+1})-\Phi^*
\le\frac1{2\mu}\dist\bigl(0,\partial\Phi(x_{t+1})\bigr)^2
\le q_{\mathrm{ad}}\bigl(\Phi(x_t)-\Phi(x_{t+1})\bigr).
\]
Since
\[
\Phi(x_t)-\Phi(x_{t+1})
=\Phi(x_t)-\Phi^*
-\bigl(\Phi(x_{t+1})-\Phi^*\bigr),
\]
we obtain
\[
(1+q_{\mathrm{ad}})\bigl(\Phi(x_{t+1})-\Phi^*\bigr)
\le q_{\mathrm{ad}}\bigl(\Phi(x_t)-\Phi^*\bigr),
\]
and therefore
\[
\Phi(x_{t+1})-\Phi^*
\le\frac{q_{\mathrm{ad}}}{1+q_{\mathrm{ad}}}
\bigl(\Phi(x_t)-\Phi^*\bigr)
=\Bigl(1-\frac1{1+q_{\mathrm{ad}}}\Bigr)
\bigl(\Phi(x_t)-\Phi^*\bigr).
\]
This is the claimed linear contraction.
\end{proof}

\subsection{Relaxed linear convergence for proximal Adam}\label{sec:pl-relaxed}

The alignment-free analysis of Section~\ref{sec:proxadam-relaxed} also admits a
linear convergence statement under the proximal PL condition. Compared with
Theorem~\ref{thm:prox-adam-pl}, the following result does not use
Assumptions~\ref{ass:align} and \ref{ass:coherence}; instead, the momentum memory is
absorbed into a Lyapunov function.

\begin{thm}[Relaxed linear convergence of proximal Adam under PL]\label{thm:prox-adam-pl-relaxed}
Let Assumptions~\ref{ass:smooth} and \ref{ass:bounded} hold, suppose that $\Phi$ is bounded below, and assume that
$\Phi$ satisfies Definition~\ref{def:pl} on the set of iterates with constant
$\mu>0$. Define
\[
a=\frac M\alpha+L,
\qquad
R_\beta>\max\left\{1,\frac{2L\beta_1}{\mu(1-\beta_1)}\right\},
\qquad
W=\frac{L R_\beta}{2(1-\beta_1)},
\]
and assume
\[
0<\alpha<\frac{2\underline m(1-\beta_1)}{L(1+\beta_1R_\beta)}.
\]
Set
\[
\widehat A=\frac{\underline m}\alpha-\frac{L}{2(1-\beta_1)},\qquad
c=\widehat A-W\beta_1,
\]
\[
K=\frac{a^2}{\mu}+\beta_1W,\qquad
R=\beta_1+\frac{2L\beta_1}{\mu R_\beta},
\]
and
\[
\rho=\frac{K/c+R}{1+K/c}.
\]
Then $0<\rho<1$ and
\[
\Phi(x_t)-\Phi^*
\le
\rho^t\bigl(\Phi(x_0)-\Phi^*\bigr).
\]
\end{thm}

\begin{proof}
Let
\[
Q_t=\|d_t\|^2,\qquad
S_t=\sum_{k=0}^{t-1}\beta_1^{t-k}Q_k,\qquad
\Xi_t=\Phi(x_t)-\Phi^*.
\]
Then
\[
S_{t+1}=\beta_1S_t+\beta_1Q_t,\qquad S_0=0.
\]

\smallskip

\noindent\textbf{Step 1: Lyapunov descent.}
From \eqref{eq:adam-relaxed-descent} and $D_t\succeq \underline mI$,
\[
\Xi_{t+1}
\le
\Xi_t-\widehat AQ_t+\frac L2S_t,
\]
where
\[
\widehat A=\frac{\underline m}\alpha-\frac{L}{2(1-\beta_1)}.
\]
Define
\[
P_t=\Xi_t+WS_t.
\]
Using $S_{t+1}=\beta_1S_t+\beta_1Q_t$,
\[
\begin{aligned}
P_{t+1}
&=\Xi_{t+1}+WS_{t+1}\\
&\le
\Xi_t-\widehat AQ_t+\frac L2S_t
+W(\beta_1S_t+\beta_1Q_t)\\
&=
\Xi_t+\left(\frac L2+\beta_1W\right)S_t
+(\beta_1W-\widehat A)Q_t.
\end{aligned}
\]
Because
\[
W\ge\frac{L}{2(1-\beta_1)},
\]
we have
\[
\frac L2+\beta_1W\le W.
\]
Also, by the condition on $\alpha$,
\[
c=\widehat A-W\beta_1>0.
\]
Hence
\[
P_{t+1}\le P_t-cQ_t.
\]
In particular, $P_t$ is nonincreasing and
\begin{equation}\label{eq:adam-relaxed-qbound}
Q_t\le\frac{P_t-P_{t+1}}c.
\end{equation}

\smallskip

\noindent\textbf{Step 2: PL residual bound.}
As in the proof of Theorem~\ref{thm:prox-adam-relaxed},
\[
\dist(0,\partial\Phi(x_{t+1}))^2
\le
2a^2Q_t
+
2L^2\frac{\beta_1}{1-\beta_1}S_t.
\]
By Definition~\ref{def:pl},
\[
\Xi_{t+1}
\le
\frac1{2\mu}\dist(0,\partial\Phi(x_{t+1}))^2
\le
\frac{a^2}{\mu}Q_t
+
\frac{L^2\beta_1}{\mu(1-\beta_1)}S_t.
\]
Write
\[
c_1=\frac{a^2}{\mu},\qquad
c_2=\frac{L^2\beta_1}{\mu(1-\beta_1)}.
\]

\smallskip

\noindent\textbf{Step 3: contraction.}
Using $S_{t+1}=\beta_1S_t+\beta_1Q_t$,
\[
\begin{aligned}
P_{t+1}
&=\Xi_{t+1}+WS_{t+1}\\
&\le
c_1Q_t+c_2S_t
+W(\beta_1S_t+\beta_1Q_t)\\
&=
(c_1+\beta_1W)Q_t
+(c_2+\beta_1W)S_t.
\end{aligned}
\]
Let $K=c_1+\beta_1W$. Since
\[
P_t=\Xi_t+WS_t\ge WS_t,
\]
we have
\[
S_t\le\frac{P_t}{W}.
\]
Combining this with \eqref{eq:adam-relaxed-qbound},
\[
P_{t+1}
\le
\frac Kc(P_t-P_{t+1})
+
\left(\beta_1+\frac{c_2}{W}\right)P_t.
\]
Therefore
\[
\left(1+\frac Kc\right)P_{t+1}
\le
\left(\frac Kc+\beta_1+\frac{c_2}{W}\right)P_t.
\]
By the definitions of $W$ and $R_\beta$,
\[
\frac{c_2}{W}
=
\frac{2L\beta_1}{\mu R_\beta}
<
1-\beta_1.
\]
Thus
\[
R=\beta_1+\frac{c_2}{W}<1.
\]
Consequently,
\[
P_{t+1}\le\rho P_t,
\qquad
\rho=\frac{K/c+R}{1+K/c}<1.
\]
Finally, $\Xi_t=\Phi(x_t)-\Phi^*\le P_t$ and $P_0=\Xi_0$, so
\[
\Phi(x_t)-\Phi^*
\le
P_t
\le
\rho^t P_0
=
\rho^t\bigl(\Phi(x_0)-\Phi^*\bigr).
\]
\end{proof}

\begin{remark}
The step-size condition for Theorem~\ref{thm:prox-adam-pl-relaxed} depends on $\mu$
through $R_\beta$. If $\mu$ is small relative to $L$, then $R_\beta$ is large and the
admissible step size is correspondingly small. This is the quantitative cost of
eliminating the momentum-alignment and residual-coherence assumptions.
\end{remark}

\section{Application to TV-regularized inverse elasticity}\label{sec:application}

We now show that the proximal-Adam framework in Section~\ref{sec:prelim}
is the abstract form of the Tikhonov--total-variation inverse
linear-elasticity problem studied in the companion analysis
\cite{zhu2023,zhu2026}.  The connection is useful in both directions: the inverse
problem supplies a concrete instance in which $f$ is a nonlinear but smooth
energy-mismatch functional and $h$ is a total-variation regularizer composed with
a box projection, while the optimization theory supplies a metric-weighted proximal
iteration whose steps have closed form or admit a small inner solver.

\subsection{Forward model and inverse formulation}

Let $\Omega\subset\mathbb R^2$ be a bounded Lipschitz domain with boundary
$\partial\Omega=\overline{\Gamma_1\cup\Gamma_2}$ and
$\Gamma_1\cap\Gamma_2=\emptyset$.  For a pair of Lam\'e coefficients
$\theta=(\lambda_{\mathrm{e}},\mu_{\mathrm{e}})$ satisfying $0<\lambda_{\mathrm{e},*}\le\lambda_{\mathrm{e}}\le\lambda_{\mathrm{e}}^*$ and
$0<\mu_{\mathrm{e},*}\le\mu_{\mathrm{e}}\le\mu_{\mathrm{e}}^*$, the Cauchy stress and linearized strain tensors are
\begin{equation}\label{eq:elastic-stress}
\sigma_\theta(v)=\lambda_{\mathrm{e}}(\nabla\cdot v)I+2\mu_{\mathrm{e}}\,\bm{e}(v),\qquad
\bm{e}(v)=\tfrac12\bigl(\nabla v+(\nabla v)^\top\bigr).
\end{equation}
Given a body force $\bm f$ and a surface traction $\bm t$, the weak form of the forward
problem is
\begin{equation}\label{eq:elastic-weak}
B_\theta(u,v)
=\int_\Omega\bigl[\lambda_{\mathrm{e}}(\nabla\cdot u)(\nabla\cdot v)
+2\mu_{\mathrm{e}}\,\bm{e}(u):\bm{e}(v)\bigr]\,dx
=\int_\Omega \bm f\cdot v\,dx+\int_{\Gamma_2}\bm t\cdot v\,ds,
\end{equation}
posed in the space $V=\{v\in H^1(\Omega;\mathbb R^2):v=0\ \text{on}\ \Gamma_1\}$.
Korn's inequality and the Lax--Milgram theorem give a unique solution $u[\theta]$.

The inverse problem determines $\theta$ from a measurement $z_{\mathrm{obs}}$ of the displacement
field, with $\|z_{\mathrm{obs}}-u[\theta^\dagger]\|_{H^1}\le\delta_{\mathrm{noise}}$ for an unknown true
coefficient $\theta^\dagger$.  The admissible set carries a bounded-variation
constraint,
\begin{equation}\label{eq:elastic-admissible}
\mathcal A=\Bigl\{\theta=(\lambda_{\mathrm{e}},\mu_{\mathrm{e}}):\lambda_{\mathrm{e},*}\le\lambda_{\mathrm{e}}\le\lambda_{\mathrm{e}}^*,\ \mu_{\mathrm{e},*}\le\mu_{\mathrm{e}}\le\mu_{\mathrm{e}}^*,\
\operatorname{TV}(\lambda_{\mathrm{e}})<\infty,\ \operatorname{TV}(\mu_{\mathrm{e}})<\infty\Bigr\},
\end{equation}
where the total variation is
\begin{equation}\label{eq:elastic-TV}
\operatorname{TV}(v)=\sup\Bigl\{\int_\Omega v\,\nabla\cdot\varphi\,dx:
\varphi\in C_c^1(\Omega;\mathbb R^2),\ |\varphi|\le1\Bigr\}.
\end{equation}
Because the forward map $\theta\mapsto u[\theta]$ is nonlinear and smoothing, the
inverse problem is ill posed; see \cite{zhu2023} for a non-uniqueness
discussion.  We use the energy-norm data mismatch
\begin{equation}\label{eq:elastic-J0}
J_0(\theta;z_{\mathrm{obs}})=\frac12 B_\theta\bigl(u[\theta]-z_{\mathrm{obs}},u[\theta]-z_{\mathrm{obs}}\bigr),
\end{equation}
and the Tikhonov--total-variation objective
\begin{equation}\label{eq:elastic-J}
J_\eta(\theta;z_{\mathrm{obs}})=J_0(\theta;z_{\mathrm{obs}})+\eta\operatorname{TV}(\theta),
\qquad
\operatorname{TV}(\theta)=\operatorname{TV}(\lambda_{\mathrm{e}})+\operatorname{TV}(\mu_{\mathrm{e}}),
\end{equation}
with $\eta>0$.

\subsection{Discrete parametrization and the composite objective}

Following the composite-material setting of \cite{zhu2023}, the coefficients
are piecewise constant on angular sectors.  With
$\ell=(\lambda_{\mathrm{e},1},\dots,\lambda_{\mathrm{e},n},\mu_{\mathrm{e},1},\dots,\mu_{\mathrm{e},n})\in\mathbb R^{2n}$, the
admissible set becomes the compact box
\begin{equation}\label{eq:elastic-box}
\mathcal C=\Bigl\{\ell\in\mathbb R^{2n}:\lambda_{\mathrm{e},*}\le\lambda_{\mathrm{e},i}\le\lambda_{\mathrm{e}}^*,\ \mu_{\mathrm{e},*}\le\mu_{\mathrm{e},i}\le\mu_{\mathrm{e}}^*\Bigr\}.
\end{equation}
The smoothed discrete total variation is
\begin{equation}\label{eq:elastic-TVeps}
\operatorname{TV}_{\varepsilon_{\mathrm{TV}}}(\ell)
=\sum_{i=1}^n\tilde r(\phi_i)\Bigl(
\sqrt{(\lambda_{\mathrm{e},i+1}-\lambda_{\mathrm{e},i})^2+\varepsilon_{\mathrm{TV}}^2}
+\sqrt{(\mu_{\mathrm{e},i+1}-\mu_{\mathrm{e},i})^2+\varepsilon_{\mathrm{TV}}^2}\Bigr),
\end{equation}
with periodic indices $\lambda_{\mathrm{e},n+1}=\lambda_{\mathrm{e},1}$ and $\mu_{\mathrm{e},n+1}=\mu_{\mathrm{e},1}$.  Let $\phi_i$ denote the angular midpoint of the $i$-th sector and let $\tilde r(\phi_i)>0$ denote the radial quadrature weight of that sector.
The discrete
regularized objective is
\begin{equation}\label{eq:elastic-Jeps}
J_{\eta,\varepsilon_{\mathrm{TV}}}(\ell)
=J_0(\theta_h(\ell);z_{\mathrm{obs}})+\eta\operatorname{TV}_{\varepsilon_{\mathrm{TV}}}(\ell).
\end{equation}
The map $\ell\mapsto\theta_h(\ell)$ is the piecewise-constant coefficient assembly: on sector $i$ it takes the value $(\lambda_{\mathrm{e},i},\mu_{\mathrm{e},i})$.

This is exactly a composite problem of the form \eqref{eq:composite}.  We set
\begin{equation}\label{eq:elastic-mapping}
f(\ell)=J_0(\theta_h(\ell);z_{\mathrm{obs}}),\qquad
h(\ell)=\eta\operatorname{TV}_{\varepsilon_{\mathrm{TV}}}(\ell)+\iota_{\mathcal C}(\ell),
\qquad
\Phi(\ell)=f(\ell)+h(\ell).
\end{equation}
For the exact TV limit $\varepsilon_{\mathrm{TV}}=0$, $\operatorname{TV}_{\varepsilon_{\mathrm{TV}}}$ is replaced by
the nonsmooth sector-to-sector variation, so $h$ is convex, proper, and lower
semicontinuous; for $\varepsilon_{\mathrm{TV}}>0$ the TV part is smooth but nonquadratic, and the
nonsmoothness is carried by the box indicator $\iota_{\mathcal C}$.  The
$L$-smoothness of $f$ on a neighborhood of $\mathcal C$ and the boundedness of its
gradients follow from the compactness of $\mathcal C$ and the differentiability
results in \cite{zhu2023}. These verify
Assumptions~\ref{ass:smooth} and \ref{ass:gradbound}; together with $\varepsilon$ in \eqref{eq:elastic-iter}, they give Assumption~\ref{ass:bounded} with $\underline m=\varepsilon$ and $M=G_\infty+\varepsilon$.

\subsection{Adjoint-free gradient and proximal iteration}

A key advantage of the energy mismatch \eqref{eq:elastic-J0} is that its gradient
does not require an adjoint solve.  For $u=u[\theta]$ and $z_{\mathrm{obs}}$ fixed,
\begin{equation}\label{eq:elastic-gradJ0}
\frac{\delta J_0}{\delta\lambda_{\mathrm{e}}}
=-\frac12\bigl(\nabla\cdot(u+z_{\mathrm{obs}})\bigr)\bigl(\nabla\cdot(u-z_{\mathrm{obs}})\bigr),\qquad
\frac{\delta J_0}{\delta\mu_{\mathrm{e}}}
=-\bm{e}(u+z_{\mathrm{obs}}):\bm{e}(u-z_{\mathrm{obs}}).
\end{equation}
After integration over each sector, \eqref{eq:elastic-gradJ0} yields
$\nabla_\ell J_0$, which is the smooth term $f$ in \eqref{eq:elastic-mapping}.  The
smoothed TV term is handled inside the weighted proximal map of $h$, so its
gradient does not enter the Adam momentum.
The proximal iteration in \eqref{eq:prox-adam} therefore reads
\begin{equation}\label{eq:elastic-iter}
\ell_{t+1}
=\prox_{\alpha h}^{D_t}\Bigl(\ell_t-\alpha D_t^{-1}\hat m_{t+1}\Bigr),
\qquad
D_t=\operatorname{diag}\bigl(\sqrt{\hat v_{t+1}}+\varepsilon\bigr).
\end{equation}
Here $\varepsilon>0$ is the same metric regularization parameter as in \eqref{eq:prox-adam}, while $\varepsilon_{\mathrm{TV}}>0$ is the TV smoothing parameter in \eqref{eq:elastic-TVeps}.
The proximal map \eqref{eq:weighted-prox} with $h$ from \eqref{eq:elastic-mapping}
is a box-constrained weighted TV-denoising subproblem.  In the diagonal metric it
is not generally a single coordinatewise threshold, but it decomposes into
(i) the metric-weighted box projection, which for a diagonal $D_t$ is the
componentwise clamping operator in Section~\ref{sec:prelim}, and
(ii) an anisotropic finite-difference TV proximal step, whose dual is a projection
onto $[-1,1]$.  A small inner Chambolle--Pock loop realizes the exact proximal map~\cite{chambolle2011},
so each outer Adam step remains a legitimate instance of
\eqref{eq:prox-adam}. 

\section{Numerical illustration}\label{sec:numerics}
\subsection{Diagonal Lasso problem}\label{sec:ex-lasso}
We test Theorems~\ref{thm:prox-norm-rate}--\ref{thm:prox-adam-pl-relaxed} on a diagonal
Lasso problem,
\[
\Phi(x)=\frac12\|Ax-b\|^2+\lambda\|x\|_1,
\qquad
A=\operatorname{diag}(s_1,\dots,s_d),
\]
with $d=40$, $s_i=0.5+\frac{i-1}{39}$, $\lambda=0.05$, and deterministic data
$b=A\mathbf 1$, where $\mathbf 1=(1,\dots,1)^\top$.  Since $A$ is diagonal, the exact
minimizer and the exact optimal value $\Phi^*=1.8993563061$ are available in closed
form, and the initial gap is $\Phi(x_0)-\Phi^*=9.7886778819$ for $x_0=0$.

The smooth part is $L$-Lipschitz with
\[
L=\max_i s_i^2=1,\qquad
\mu=\min_i s_i^2=0.25,
\]
so the strongly convex function $\Phi$ satisfies the proximal PL inequality of
Definition~\ref{def:pl} with constant $\mu=0.25$.  With metric regularization
$\varepsilon=1$, every $D_t$ satisfies $D_t\succeq I$, hence $\underline m=1$.  On the reported
trajectory the coordinates of the iterates remain in $[0,1]$; consequently
$\|\nabla f(x_t)\|_{\infty}\le1$, the adaptive metric satisfies
$D_t\preceq2I$, and Assumptions~\ref{ass:bounded} and~\ref{ass:gradbound} hold with
$M=2$.  All constants quoted below are the ones supplied by the corresponding theorem.

We run two algorithm configurations.  The first is proximal normalized gradient with
$\alpha=0.9$, $\beta_2=0.99$, and $\varepsilon=1$, so that
\[
\alpha=0.9\le\frac{\underline m}L=1.
\]
The second is proximal Adam with
\[
\alpha=0.05,\qquad \beta_1=0.2,\qquad \beta_2=0.99,\qquad \varepsilon=1,
\]
and $T=356$ iterations.  This configuration is deliberately chosen so that both the
aligned and the relaxed Adam analyses can be verified on the same trajectory.

For the aligned analysis we estimate the momentum-alignment and residual-coherence
constants over the reported non-asymptotic trajectory
$\{\,t:\Phi(x_t)-\Phi^*>10^{-4}\,\}$:
\[
\kappa_{\mathrm{emp}}=5.1796\times10^{-3},
\qquad
C_{0,\mathrm{emp}}=1.9790\times10^{3}.
\]
They satisfy the compatibility condition of Theorem~\ref{thm:prox-adam},
\[
\frac{\kappa\sqrt{C_0}}{\underline m}+\frac{\alpha L}{2\underline m}
=0.2554
\le\frac12.
\]
For the relaxed analysis, Theorem~\ref{thm:prox-adam-relaxed} requires
\[
\alpha\le\frac{\underline m(1-\beta_1)}{L(1+\beta_1)}
=\frac{1-\beta_1}{1+\beta_1}=0.6667,
\]
which holds.  Choosing $R_\beta=2.1$ in
Theorem~\ref{thm:prox-adam-pl-relaxed}, the admissible step-size bound is
\[
\alpha<\frac{2\underline m(1-\beta_1)}{L(1+\beta_1R_\beta)}=1.1268,
\]
which also holds.

The objective gap is fitted as
\[
\log\Delta_t\approx r_{\mathrm{fit}}\,t+c_{\mathrm{fit}},
\qquad
\rho_{\mathrm{emp}}=\exp(r_{\mathrm{fit}}),
\]
with the machine-zero plateau removed.  Table~\ref{tab:lasso-theorem} summarizes
the theoretical constants and the empirical observations.  The empirical contraction
factor for proximal normalized gradient is $0.5625$, while Theorem~\ref{thm:prox-pl}
guarantees the upper contraction factor $q/(1+q)=0.9771$.  The same proximal Adam run
gives the empirical contraction factor $0.9690$.  The aligned PL theorem gives the
contraction factor $q_{\mathrm{ad}}/(1+q_{\mathrm{ad}})=0.9976$, and the relaxed PL
theorem gives $\rho=0.9999$; both are strict contractions even though they are
conservative.  Thus the observed geometric decay is predicted by
Theorems~\ref{thm:prox-pl}, \ref{thm:prox-adam-pl}, and
\ref{thm:prox-adam-pl-relaxed}, respectively.

\begin{table}
\centering\small
\setlength{\tabcolsep}{5pt}
\renewcommand{\arraystretch}{1.35}
\begin{tabularx}{\textwidth}{@{}>{\raggedright\arraybackslash}p{3.2cm} c L c L@{}}
\toprule
Configuration & Theorems & Theoretical constants & Predicted contraction & Empirical result\\
\midrule
Normalized gradient ($\alpha=0.9$) & \ref{thm:prox-norm-rate}, \ref{thm:prox-pl} &
$C=21.4$ & $\le 0.9771$ & $2.45\times10^{-27}$; rate $0.5625$\\
\addlinespace
Aligned Adam ($\alpha=0.05$, $\beta_1=0.2$) & \ref{thm:prox-adam}, \ref{thm:prox-adam-pl} &
$C_{\mathrm{ad}}=45.6$ & $\le 0.9976$ &
$5.78\times10^{-5}$; rate $0.9690$\\
\addlinespace
Relaxed Adam ($R_\beta=2.1$) & \ref{thm:prox-adam-relaxed}, \ref{thm:prox-adam-pl-relaxed} &
$C_{\mathrm{mom}}=336.2$ & $\le 0.9999$ &
$5.78\times10^{-5}$; rate $0.9690$\\
\bottomrule
\end{tabularx}
\caption{Theorem-to-configuration verification in the deterministic diagonal Lasso
example.  The stationarity constants are those in the corresponding theorem
statements; the remaining constants are the empirical values on the reported
trajectory.}
\label{tab:lasso-theorem}
\end{table}

The stationarity checks in Table~\ref{tab:lasso-theorem} are consistent with the
explicit $O(1/N)$ bounds.  For instance, the Theorem~\ref{thm:prox-norm-rate} bound at
$N=120$ is $1.7438$, while the observed cumulative minimum residual is
$2.4511\times10^{-27}$.  For the aligned Adam run the Theorem~\ref{thm:prox-adam} bound
at $N=356$ is $5.7196$, and the relaxed Theorem~\ref{thm:prox-adam-relaxed} bound is
$9.2446$, both above the observed value $5.7822\times10^{-5}$.

\begin{figure}[H]
\centering
\includegraphics[width=0.95\textwidth]{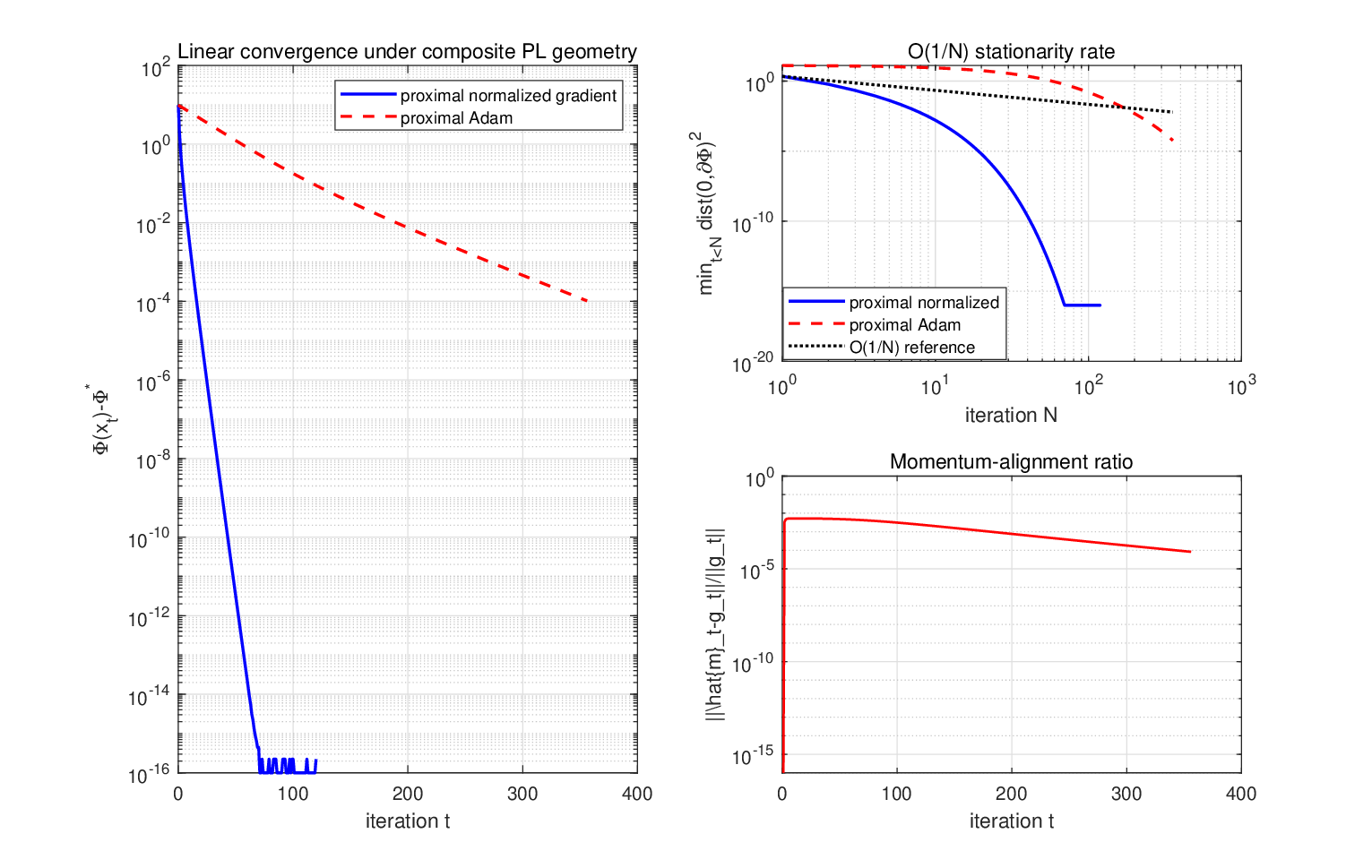}
\caption{ Left: objective gap 
$\Phi(x_t)-\Phi^\ast$ for proximal normalized gradient (solid) and proximal Adam
(dashed). Top right: cumulative minimum stationarity residual
$\min_{t<N}\dist(0,\partial\Phi(x_{t+1}))^2$ together with the $O(1/N)$ reference
line. Bottom right: the Adam momentum-alignment ratio. }
\label{fig:composite}
\end{figure}

Figure~\ref{fig:composite} collects the numerical checks.  The left panel shows
geometric decay of the objective gap.  Proximal normalized gradient reduces the gap
from $9.7887$ at $t=0$ to $4.5027\times10^{-10}$ at $t=40$, and the compatible Adam run
reduces it from $9.7887$ to $1.0270\times10^{-4}$ at $t=356$; at larger $t$ it
continues to decrease to machine precision.  The top-right panel shows that the
cumulative minimum stationarity residual follows the $O(1/N)$ reference line in
log--log scale, confirming the stationarity statements of
Theorems~\ref{thm:prox-norm-rate}--\ref{thm:prox-adam-relaxed}.  The bottom-right
panel reports the momentum-alignment ratio; its smallness is the basis of the
aligned verification in Theorems~\ref{thm:prox-adam} and
\ref{thm:prox-adam-pl}.

\subsection{Two-dimensional image TV denoising}\label{sec:ex-denoise}

The second example tests the weighted proximal operator on a genuine
nonsmooth regularizer in two dimensions.  We denoise a $64\times64$ synthetic
image $b$ with the anisotropic total-variation objective
\[
\Phi(x)=\frac12\|x-b\|^2+\lambda\operatorname{TV}(x),
\qquad
\operatorname{TV}(x)=\sum_{i,j}|x_{i+1,j}-x_{i,j}|
+\sum_{i,j}|x_{i,j+1}-x_{i,j}|,
\]
with periodic boundary conditions and $\lambda=0.2$.  The image contains a
background level, two disks, and a vertical bar, corrupted by additive Gaussian
noise of standard deviation $0.18$.  The metric-weighted proximal step
\[
\prox_{\alpha h}^{D}(z)=\arg\min_y\Bigl\{\lambda\operatorname{TV}(y)
+\frac1{2\alpha}\|y-z\|_D^2\Bigr\}
\]
is not a coordinatewise threshold. 
Both proximal normalized gradient and proximal Adam are run with $\alpha=1$,
$\beta_2=0.99$, $\varepsilon=1$, and $\beta_1=0.5$ for Adam.

The weighted proximal step is solved by a Chambolle--Pock inner loop with the
correct componentwise weighted update
\[
x\leftarrow
\frac{x-\tau_{\mathrm{in}}K^\top p+\frac{\tau_{\mathrm{in}}}{\alpha}D\odot z}
{1+\frac{\tau_{\mathrm{in}}}{\alpha}D},
\qquad
p\leftarrow\operatorname{clip}\bigl(p+\sigma K\bar x,[-\lambda,\lambda]\bigr),
\]
where $K$ denotes the discrete gradient operator.  The dual variables are
warm-started from the previous outer iteration and the loop is terminated when
the dual fixed-point residual is below $10^{-8}$.  With this inner solve the
proximal normalized-gradient objective is nonincreasing up to the inner-solver
tolerance, as required by the composite descent lemma.  Both methods start from
$x_0=0$, which is distinct from the observation $b$, and the stationarity and
cumulative-minimum residuals follow the $O(1/N)$ reference line, consistent
with Theorems~\ref{thm:prox-norm-rate} and \ref{thm:prox-adam}.
The noisy observation has PSNR $14.9232$ dB, while the two proximal methods
attain $30.0374$ dB, with final objectives $7.362560\times10^1$ and
$7.362558\times10^1$ and minimum stationarity residuals $1.586\times10^{-7}$
and $2.348\times10^{-8}$, respectively.  The equality of the two PSNR values is
expected because the composite objective is strongly convex and both runs converge
to the same TV minimizer.
\IfFileExists{tv_denoising_result.eps}{%
\begin{figure}[H]
\centering
\includegraphics[width=0.95\textwidth]{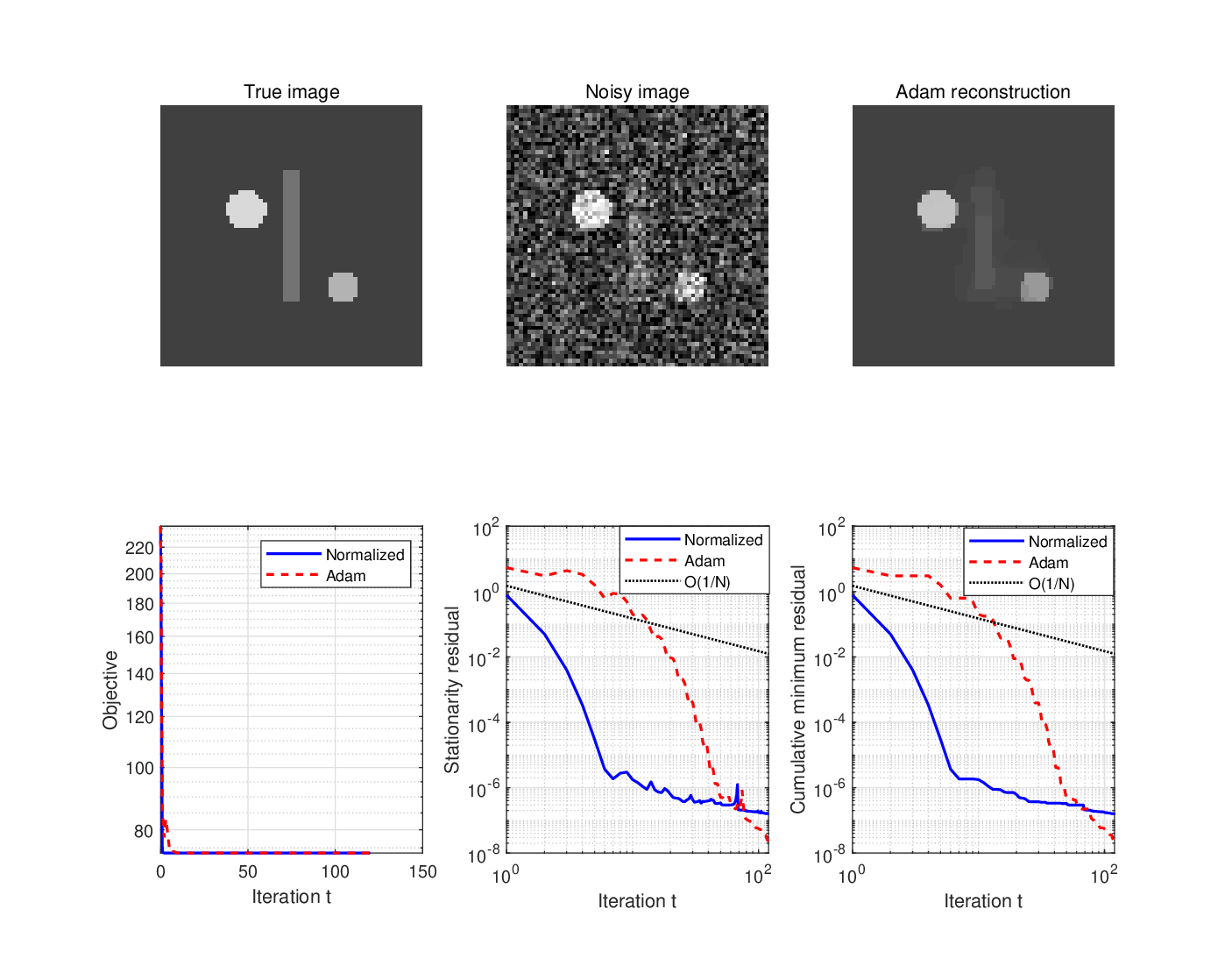}
\caption{Two-dimensional image TV denoising with the corrected weighted
Chambolle--Pock inner solver.  Top row: true image, noisy observation, and
proximal Adam reconstruction.  Bottom row: objective values, stationarity
residual, and cumulative minimum residual with the $O(1/N)$ reference.}
\label{fig:tv-denoise}
\end{figure}
}{}

\subsection{Two-dimensional elastic parameter identification}\label{sec:ex-elastic}

The third example is the inverse-elasticity application of
Section~\ref{sec:application}.  On the unit square with homogeneous Dirichlet
boundary conditions we solve the plane-strain forward problem
\[
-\nabla\cdot\sigma(u)=\bm f,\qquad
\sigma(u)=\lambda_{\mathrm{e}}(\nabla\cdot u)I+2\mu_{\mathrm{e}}\,\bm{e}(u),
\]
by a second-order finite-difference discretization on a $24\times24$ interior
grid.  Two body forces $\bm f^{(1)},\bm f^{(2)}$ generate two noisy displacement
measurements, and the unknown parameter vector $\ell=(\lambda_{\mathrm{e}},\mu_{\mathrm{e}})$ is recovered
from the normalized least-squares mismatch
\[
f(\ell)=\frac12\Bigl(
\frac{\|u^{(1)}(\ell)-z_{\mathrm{obs}}^{(1)}\|^2}{\|z_{\mathrm{obs}}^{(1)}\|^2}
+\frac{\|u^{(2)}(\ell)-z_{\mathrm{obs}}^{(2)}\|^2}{\|z_{\mathrm{obs}}^{(2)}\|^2}\Bigr),
\qquad
h(\ell)=\iota_{[0.2,5]^2}(\ell).
\]
Because the metric $D_t$ is diagonal and $h$ is a box indicator, the weighted
proximal map \eqref{eq:weighted-prox} is the componentwise clamping operator
in Section~\ref{sec:prelim}.  We start from
$\ell_0=(1.0,0.5)$ with true coefficients $(\lambda_{\mathrm{e}},\mu_{\mathrm{e}})=(1.2,0.8)$.

\IfFileExists{elastic_inverse_result.eps}{%
\begin{figure}[H]
\centering
\includegraphics[width=0.95\textwidth]{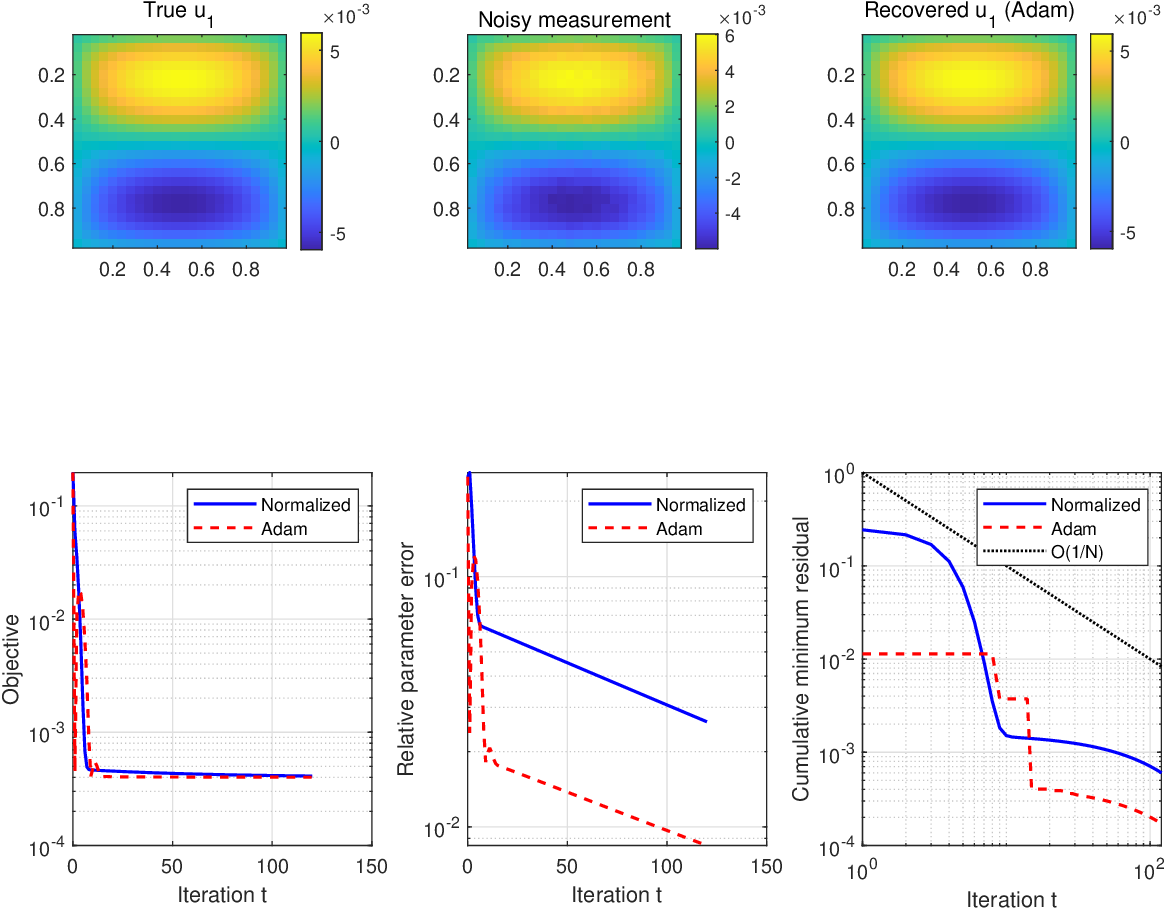}
\caption{Two-dimensional elastic parameter identification.  Top row: true,
noisy, and recovered first displacement component.  Bottom row: objective decay,
relative parameter error, and cumulative minimum stationarity residual.}
\label{fig:elastic-inverse}
\end{figure}
}{}

\paragraph{Interpretation}
With $2\%$ relative noise, proximal normalized gradient recovers
$(\lambda_{\mathrm{e}},\mu_{\mathrm{e}})=(1.1642,0.8126)$ and proximal Adam recovers
$(1.1887,0.8046)$.  The Adam relative parameter error $8.429\times10^{-3}$ is
smaller than the normalized-gradient error $2.630\times10^{-2}$.  In this example
the Adam momentum smooths the adaptively normalized direction, whereas the
normalized-gradient method uses only the instantaneous normalized gradient.
The final mismatch values $4.1069\times10^{-4}$ and $4.0072\times10^{-4}$ are
both close to the noise floor
$\frac12\bigl(0.02^2+0.02^2\bigr)=4\times10^{-4}$, showing that the two
schemes fit the data without overfitting.  The minimum stationarity residuals
are $6.020\times10^{-4}$ and $1.716\times10^{-4}$, confirming convergence of
the projected adaptive iteration.
\section{Conclusion}\label{sec:conclusion}
This paper developed a nonasymptotic convergence framework for proximal
normalized gradient and proximal Adam applied to the composite problem
\begin{equation}
\min_x\ \Phi(x)=f(x)+h(x),
\end{equation}
where $h$ is a convex prox-friendly regularizer and $f$ is a possibly
nonconvex smooth term with $L$-Lipschitz gradient.  The analysis rests on three core tools: a metric-weighted composite descent lemma, a
residual--step equivalence that converts weighted step decay into
subgradient stationarity, and the proximal PL condition introduced in
Definition~\ref{def:pl}; the relaxed Adam results additionally use the momentum-error estimate of Lemma~\ref{lem:adam-momentum-error}.  Under bounded positive-definite metrics and
boundedness of $\Phi$ from below, Theorems~\ref{thm:prox-norm-rate}
and~\ref{thm:prox-adam} show that the minimum subdifferential-distance
residual decays as $O(1/N)$; the Adam version uses momentum-alignment and
residual-coherence assumptions together with the metric-compatibility
condition $\frac{\kappa\sqrt{C_0}}{\underline m}+\frac{\alpha L}{2\underline m}\le\frac12$ to control the bias of the adaptive
direction. Theorem~\ref{thm:prox-adam-relaxed} removes these two assumptions by
estimating the momentum error directly, at the cost of a more restrictive step-size
range and a larger constant.
Theorems~\ref{thm:prox-pl} and~\ref{thm:prox-adam-pl} upgrade the results
to linear convergence with explicit contraction factors
$1-1/(1+q)$ and $1-1/(1+q_{\mathrm{ad}})$, respectively. The relaxed
Theorem~\ref{thm:prox-adam-pl-relaxed} gives the corresponding linear rate
without momentum-alignment and residual-coherence, using a Lyapunov function that
absorbs the momentum memory.

Three numerical examples support the theory.  The diagonal Lasso test
confirms the $O(1/N)$ stationarity bound and the proximal PL linear
contraction for proximal normalized gradient, with empirical linear rate
$0.5625$; the compatible proximal Adam run has empirical rate $0.9690$ and
satisfies the metric-compatibility condition of
Theorem~\ref{thm:prox-adam}, while the relaxed analysis verifies the same run
through Theorem~\ref{thm:prox-adam-pl-relaxed}. Two-dimensional image TV
denoising verifies the weighted proximal operator with a corrected
Chambolle--Pock inner solver and produces a nonincreasing objective up to the
inner-solver tolerance.  In the
two-dimensional elastic parameter identification, with true coefficients
$(\lambda_{\mathrm{e}},\mu_{\mathrm{e}})=(1.2,0.8)$, proximal normalized gradient recovers
$(1.1642,0.8126)$ and proximal Adam recovers $(1.1887,0.8046)$,
corresponding to relative parameter errors $2.630\times10^{-2}$ and
$8.429\times10^{-3}$.  The final data mismatches $4.1069\times10^{-4}$ and
$4.0072\times10^{-4}$ lie close to the noise floor, indicating that neither
scheme overfits the data; the Adam smaller parameter error reflects momentum smoothing
on top of the coordinatewise normalization used by both schemes.

Natural extensions are 
accelerated proximal-Adam schedules, sharper momentum-alignment and
residual-coherence constants, nonconvex regularizers analyzed through
KL-type or desingularizing geometry.

\section*{Acknowledgments}

This work was partially supported by the NSFC Projects No. 12526566 and Liaoning Provincial Natural Science Foundation Project No. 2025-BS-0410.

\section*{Declarations}

\noindent\textbf{Funding.} This work was partially supported by the NSFC Projects No. 12526566 and Liaoning Provincial Natural Science Foundation Project No. 2025-BS-0410.

\noindent\textbf{Data Availability.} This paper has no associated data.

\noindent\textbf{Competing Interests.} The authors have no relevant financial or non-financial interests to disclose.

\noindent\textbf{Author Contributions.} All authors contributed to the study conception and design. Theoretical analysis and numerical experiments were performed by Xiaopeng Zhu. The first draft of the manuscript was written by Xiaopeng Zhu; Zhongyi Huang reviewed and revised the manuscript. All authors read and approved the final manuscript.

\end{document}